\documentclass{amsart}
\usepackage{amsmath,amsfonts,amssymb,psfig}

\newtheorem{theorem}{Theorem}[section]

\newtheorem{lemma}[theorem]{Lemma}
\newtheorem{corollary}[theorem]{Corollary}

\newtheorem{conjecture}[theorem]{Conjecture}

\theoremstyle{definition}
\newtheorem{remark}[theorem]{Remark}
\newtheorem{example}[theorem]{Example}

\newcommand{\abs}[1]{{\left|#1\right|}}
\newcommand{\cal}{\mathcal}

\newcommand{\Z}{{\mathbf{Z}}}
\newcommand{\N}{{\mathbf{N}}}
\newcommand{\Q}{{\mathbf{Q}}}

\newcommand{\Zp}{{\Z_p}}

\newcommand{\nm}{^{\mathrm{nm}}}
\newcommand{\bt}{\beta}
\newcommand{\tf}{\tilde{f}}
\newcommand{\ov}{\overline}
\newcommand{\PP}{{\mathbf{P}}}
\newcommand{\pr}{\mathrm{pr}}
\newcommand{\spc}{\mathrm{\ splits\ completely}}
\newcommand{\mnc}{\rm\ is\ monic}
\newcommand{\da}{d_{\alpha}}
\newcommand\floor[1]{\lfloor#1\rfloor}
\newcommand{\TT}{\gamma}
\newcommand{\TTT}{\nu}
\newcommand{\FF}{F}
\newcommand{\GG}{G}

\title{The probability that a random monic $p$-adic 
polynomial splits into
linear factors}

\date{\today}

\begin{document}
\author{
Joe Buhler \and
Daniel Goldstein \and 
David Moews \and 
Joel Rosenberg}
\address
{Center for Communications Research\\
4320 Westerra Court\\
San Diego, CA 92121}
\email{\tt \{buhler,dgoldste,dmoews,joelr\}@ccrwest.org}
\maketitle
\tableofcontents

\begin{abstract}
Let $n$ be a positive integer and let $p$ be a prime.
We calculate the probability that a random 
monic polynomial 
of degree $n$ with coefficients in the ring $\Zp$ of $p$-adic
integers splits over $\Zp$ into linear factors.
\end{abstract}

\section{Introduction}\label{intro}
Let $R$ be a complete discrete valuation ring with finite residue field.
The main result of this note, Theorem~\ref{thm:1},
is a formula, for each 
positive integer $n$, for the
probability that a random monic polynomial
of degree $n$ with coefficients in $R$ splits over $R$ into linear factors.

By a standard result, $R$ is either 
the ring of formal power series over a finite field, 
or the ring of integers of a finite extension of the field of 
$p$-adic numbers for some prime $p$.
In either case we write $k$ for the residue field of $R$ and set $q=|k|$.
We encourage the reader to focus on the 
special cases $R=\Z_p$ and $GF(p)[[t]]$, as 
the statements and proofs in the general case are essentially the same
as in these two special cases.

Given  $(c_1,\dots,c_n)$ in $R^n$, we form the monic
(i.e. leading coefficient~1) polynomial
$f(x)= x^n + c_1x^{n-1} + \dots+c_n$.
Let $P_n(R)$ be the set of monic degree $n$ 
polynomials with coefficients in $R$,
and identify $P_n(R)$ with $R^n$ as above.
Let $S_n(R)\subseteq P_n(R)$ consist of those polynomials
that split over $R$ into linear factors (or equivalently, 
have $n$ roots in $R$, counting multiplicity).  We identify
$S_n(R)$ with a subset of $R^n$.
In Section~\ref{secff}
we will similarly write $P_n(k)$ (resp.\ $S_n(k)$)
for the set of monic polynomials of degree $n$ (resp.\ monic polynomials
of degree $n$ that split completely) over the residue field $k$.

The compact abelian group $R$ has a Haar measure $\mu$ of total mass~1.
Give $R^n$ the product topology and product measure $\mu_n$
(alternatively, view $R^n$ as a compact abelian group in its own right).
One can show that $S_n(R)$ is $\mu_n$-measurable as follows.
Given $(a_1,\dots,a_n)\in R^n$, we can form the polynomial 
$\prod_{1\le i\le n}(x-a_i)$. The resulting map $\Xi$ from 
$R^n$ to $P_n(R)$ is continuous and has image $S_n(R)$.
Thus, $S_n(R)$ is compact, being the continuous image of a compact set,
whence closed and in particular measurable.

A main result of this note is a recursive formula for
$r_n:=\mu_n(S_n(R)).$
Loosely speaking, this is the probability that a random monic polynomial
of degree $n$ over $R$ splits completely over $R$.
Let $q$ be the order of the residue field of $R$.

The first four values of $r_n$ are:
\begin{eqnarray*}
 r_1 &=&1,\\
 r_2 &=&q/{2(q+1)},\\
 r_3 &=&(q^2-q+1)(q-1)q^3/ {6(q+1) (q^5-1)},\\
 r_4 &=&h(q)(q-1)^4q^6/{24 (q^2-1)^2 (q^5-1) (q^9-1)},
\end{eqnarray*}
where $h(q)= q^8 - 2q^7 + q^6 + 2q^5 - q^4 + 2q^3 + q^2 - 2q + 1$.

We also give a formula for the probability that a random polynomial that
is not necessarily monic of degree $n$ splits over $R$.

A student of the first-named author, Asher Auel \cite{Auel}, 
working independently, 
has calculated the Jacobian determinant of the map $\Xi$,
then calculated the resulting $n$-fold $p$-adic integral.
His approach offers an alternative to the method of this note.

Here is our main result.
Set $r_0=1$. 
For a sequence  $d=(d_0,d_1,\dots,d_{q-1})$ of nonnegative integers, write 
$\abs{d}$ for their sum, $d_0+d_1+\cdots+d_{q-1}$.
\begin{theorem}\label{thm:1}
Let $R$ be a complete discrete valuation ring with finite residue
field $k$.  Set $q=\abs{k}$. For $n$ a positive integer, write $r_n$ for
the probability that a random degree $n$ monic polynomial over $R$
splits into linear factors, and set $r_0=1$.  If $n\ge0$, then
\begin{equation}
\label{maineq}
r_n=\sum_{\abs{d}=n}\ \prod_{0\le i\le q-1} q^{-\binom{d_i+1}{2}}r_{d_i}.
\end{equation}
\end{theorem}

\begin{remark}
Note that $r_n$ appears in 
$q$ terms on the right,
each time with coefficient $q^{-\binom{n+1}2}.$ 
But $0=1-\binom{n+1}{2}$ if and only if $n=1$.
So, as long as $n>1$, 
we can solve for $r_n$
in terms of $r_0,\dots,r_{n-1}$.
\end{remark}

Here is an approach to $r_n$ using power series.
We define $\FF(t)$ and $\GG(t)$ in $\Q[[t]]$ by the
formulas
\begin{eqnarray}\label{eq1}
\FF(t) &= &\sum_{n\ge0}r_n t^n, \mbox{ and}\\\label{eq2}
\GG(t) &= &\sum_{n\ge0}\frac{r_n}{q^{\binom{n +1}{2}}}t^n.
\end{eqnarray}
Then Theorem~\ref{thm:1} gives that 
\begin{equation}
\label{eq3}
\FF = \GG^q.
\end{equation}

We use a method due to Euler to 
get an efficient method for computing $r_n$.

\begin{corollary}\label{cor:1}
For $n\ge0$, we have
$$
\sum_{0\le j\le n}(n-(q+1)j)r_{n-j} \frac{r_j}{q^{\binom{j+1}{2}}}=0.
$$
In particular, $r_n$ is a rational function in $q$. 
\end{corollary}

\section{Proofs}
For $r\in R$, consider the map $\psi_r$ that takes the monic degree $n$ 
polynomial $f(x)$ to
$f(x-r)$. Thus, $\psi_r$ maps $P_{n}(R)$ to itself.

\begin{lemma}\label{lem:0}Let $r\in R$.
The map $\psi_r$ from $P_{n}(R)$
to itself preserves the measure $\mu_n$.
\end{lemma}
\begin{proof}
Our identification of $P_n(R)$ with $R^n$
makes $\psi_r$ into an affine map.
For, write $f=x^n + f_1$ with $\deg(f_1)<n$.
Then $\psi_r$ is the composition of the linear map
$\psi'_r(x^n + f_1)= x^n + f_1(x-r)$ and the translation map
$f\mapsto f+ (x-r)^n -x^n$.

Since any translation preserves measure, it is enough to check the
linear map $\psi'_r$.

A basis for $P_n(R)$ is given by the set of 
monomials $x^i$ for $0\le i\le n-1$.
In terms of this basis, the matrix of $\psi'_r$ is upper triangular with ones 
on the diagonal.  Hence it has Jacobian determinant equal to 1.
\end{proof}

We proceed by a series of lemmas. The basic idea is as follows.
Suppose the monic polynomial $f(x)$ in $P_n(R)$
splits completely over $R$. 
Write $\pi R$ for the maximal ideal of $R$, and write
$\ov{f}$ for the reduction of $f$ mod $\pi R$.
Then $\ov{f}$ splits completely over $k=R/\pi R$.
Write this factorization as 
\begin{equation}
\ov{f}(x) = \rho(x) = \prod_{\alpha\in k} \rho_\alpha(x),
\qquad \rho_\alpha(x) = (x-\alpha)^{\da}.
\end{equation}

The polynomials $\rho_\alpha(x) = (x-\alpha)^{\da}$ are pairwise relatively
prime.
Hence by Hensel's lemma and induction, $f=\prod_{\alpha\in k}f_{\alpha}$, 
where $f_{\alpha}$ in $R[x]$ reduces to $\rho_\alpha$ in $k[x]$.

Of course, the probability that $f$ reduces to 
$\rho$
(or to any given degree $n$ monic polynomial in $k[x]$)
is exactly $q^{-n}$.
The desired result will follow from the following lemma.
\begin{lemma}\label{lem:1} \
\begin{enumerate}
\item $\pr(f \spc|\ov{f}=\rho)=\prod_{\alpha\in k}
\pr(f_{\alpha}\spc|\ov{f_{\alpha}}=\rho_\alpha)$.
\item $\pr(f_{\alpha}\spc|\ov{f_{\alpha}}=\rho_\alpha)=
q^{-\binom{\da}{2}}r_{\da}.$
\end{enumerate}
\end{lemma}

To prove Theorem \ref{thm:1} from the lemma,
it suffices to choose a bijection from $k$ to the set $S=\{0,1,\dots,q-1\}$
and to use that $d_i = \binom{d_i+1}{2}-\binom{d_i}{2}$
and that $\abs{d}=\sum_{i\in S}d_i$.

Let $g$ be a polynomial in $P_n(k)$.
Write $P_g$ for the set of polynomials $f$ in $P_n(R)$ whose reduction
$\ov{f}$ is $g$.  As $P_g$ is a measurable subset of $P_n(R)$, it inherits a 
topology from $P_n(R)$, and also a measure $\mu_n$, which we renormalize so 
that $P_g$ has total measure~1.  Part~(1) of Lemma~\ref{lem:1} now follows from:

\begin{lemma}\label{lem:2}(Measure theoretic version of Hensel's lemma).
Let $g$ and $h$ in $k[x]$ be relatively prime monic polynomials of degree
$m$ and $n$, respectively.  The multiplication map
$$
M: P_g\times P_h\to P_{gh}
$$
is a bijection, in fact an isomorphism of topological spaces, by the 
usual Hensel's lemma.
Then $M$ preserves measure, that is, $M_{\ast}(\mu_m\times\mu_n)=\mu_{m+n}$.
\end{lemma}

The lemma is an immediate consequence of the following, somewhat more
general result.

\begin{lemma}\label{lem:3}
Let 
$A=\mbox{}\prod_{1\le i}a_i$ and
$B=\mbox{}\prod_{1\le i}b_i$ be countable products of finite sets.
Normalize counting measure so that
each of $a_i$ and $b_i$ has total mass~1, for all $i$, 
and give $A$ and $B$ the product measure.
Suppose there is a compatible system of bijections $\phi_n$
between the partial products $A_n=\mbox{}\prod_{1\le i\le n}a_i$
and  $B_n=\mbox{}\prod_{1\le i\le n}b_i$.
That is, for each $n\ge1$, the map $\phi_n$ 
from $A_n$ to $B_n$ is bijective, and 
for $m\le n$, if $y$ in $A_n$ has $x$ in $A_m$ as initial string, 
then $\phi_n(y)$ has $\phi_m(x)$ as initial string.
Then the map $\phi=\displaystyle\lim_{\leftarrow}\phi_n$ from 
$A$ to $B$ is a measure-preserving bijection.
\end{lemma}
\begin{proof}Certainly $\phi$ is bijective.
Since $\phi$ takes basic open sets to basic open sets of the
same volume, it follows that $\phi$ preserves measure.
\end{proof}

{\it Proof of part~(2) of Lemma~\ref{lem:1}}.
By Lemma~\ref{lem:0}, it suffices to treat the case $\alpha=0$.
Set $f=f_0.$ We have $f$ in $R[x]$ monic with reduction $x^n$ in $k[x]$.
Assume for the moment that $f=x^n + c_1x^{n-1} + \dots+c_n$ does split completely, say as
$$
f = (x-a_1) \cdots (x-a_n).
$$
A root of $f$ reduces to a root of $\ov{f}$, so that $a_i$ lies in $\pi R$ 
for each $i$,
whence a necessary condition for $f$ to split completely is: 
$c_i$ lies in $\pi^i R$ for each $i$.
The probability that this condition holds
is $q^{-\binom n2},$ the product of $q^{-i}$ for $0\le i\le n-1$.
Set 
$$
\tf=f(\pi x)/\pi^n.
$$

The necessary condition is met if and only if $\tf$ has coefficients in $R$.
Conditioned on its being met,
$\tf$ is distributed like a random polynomial in $P_n(R)$, and
the result follows, because 
$f$ splits completely if and only if $\tf$ splits completely.
\qed
\medskip

Corollary~\ref{cor:1} follows from the following efficient recursive
computation of the coefficients of a high power of a known power series
which is due to Euler~\cite{euler}.
The authors are grateful to H.~Wilf for pointing out Lemma~\ref{lem:w}.

\begin{lemma}\label{lem:w}
Let $L$ be a field. 	
Assume either (1) $q\in L$ and $L$ has characteristic 0, or (2) $q$ is an 
integer.
Let $\FF=\sum_{n\ge 0}\FF_nt^n$ and $\GG=\sum_{n\ge0}\GG_nt^n$ in $L[[t]]$ be 
formal power series with constant term~1
that satisfy $\FF=\GG^q.$
Then 
\begin{equation}
\label{w}
\sum_{0\le j\le n}(n-(q+1)j)\FF_{n-j}\GG_j=0.
\end{equation}
\end{lemma}
\begin{proof}
Take logarithmic derivatives to get 
$$
\frac{\FF'}{\FF}=q\frac{\GG'}{\GG}.
$$
Next, cross-multiply to get
$$
\FF' \GG t=q \FF \GG' t,
$$
and equate coefficients of $t^n$ to get
$$
\sum_{0\le j\le n-1}(n-j)\FF_{n-j}\GG_j = \sum_{1\le j\le n}q\FF_{n-j}j\GG_j.
$$
There is no harm in including the term $j=n$ in the sum on the left and
the term $j=0$ in the sum on the right, 
since these terms are zero.
Now subtract to get equation~\eqref{w}.
\end{proof}
This proves Corollary~\ref{cor:1}.
\section{Finite fields}
\label{secff}
Recall that $R$ is a complete discrete valuation ring with 
finite residue field $k$, where $q=\abs{k}$.
Let $n$ be a positive integer.
Write $\bar r_n$ for the probability that a random 
monic degree $n$ polynomial with coefficients in $k$
splits completely. Then $\bar r_n$, as well as $r_n$, depends only on 
$q$ and $n$. Our convention is that $r_0=\bar r_0=1$.
We note some properties of $r_n$ and $\bar r_n.$

1. We have $r_n\le \bar r_n$. For, if 
a monic polynomial $f$ in $R[x]$ splits completely, then
$\ov{f}$ splits completely in $k[x]$.

2. We have $\sum_{n\ge0}\bar r_nt^n=(1-t/q)^{-q}$.
For, by the result sometimes called the
stars and bars theorem,
the number of monic degree $n$ polynomials that
split completely is 
$q^n\bar r_n= \binom{n+q-1}{q-1}$.
Now use the 
binomial expansion $(1-t)^{-q}=\sum_{n\ge0}\binom{-q}{n}(-t)^n$,
and the fact that $\binom{-q}{n}=(-1)^n\binom{q-1+n}{n}$.

3. We have 
$\lim_{q\to\infty} r_n=1/{n!}=\lim_{q\to\infty} \bar r_n$.
The second equality follows from Observation~2.
For $f$ monic in $R[x]$, the probability that $\ov{f}$ has a repeated
root is at most $1/q$. (Proof: the set $S$ of monic polynomials $f$
with constant and linear term in $\pi R$ has measure $1/q^2$. Hence the
union of the sets $\psi_r(S)$ over a set of lifts $r$ 
of the elements of $k$ has
measure at most $1/q$.)
But this tends to zero as $q\to\infty$.  This proves the first equality.

\section{Non-monic polynomials}
In this section we drop the assumption that our polynomials are
monic.
We identify 
the set $P\nm_n(R)$ of degree $\le n$ polynomials
with coefficients in $R$ with $R^{n+1}$. 
It gets the measure $\mu\nm_n = \mu^{n+1}$. 
Let $S\nm_n(R)$ be those polynomials $f(x)$ in $P\nm_n(R)$ 
that can factored into $n$ linear factors, i.e. can be written in the form
\begin{equation}
\label{sp}
f = (b_1x-a_1)\cdots(b_nx-a_n),
\end{equation}
with $a_i$, $b_i\in R$.

One shows that $S\nm_n(R)$ is measurable as follows. Take a pair 
$(a_1,\dots,a_n)$ and $(b_1,\dots,b_n)$ of $n$-tuples and construct the 
polynomial $f$ as in \eqref{sp}.
The resulting continuous map from the compact set 
$R^{2n}$ to $P\nm_n(R)$ has image
$S\nm_n(R)$. Thus $S\nm_n(R)$ is compact, hence closed.

The goal of this section is a formula for $r\nm_n=\mu\nm_n(S\nm_n(R))$.
This may be interpreted as the probability that a 
random, not necessarily monic, polynomial of degree $\le n$ 
splits completely into linear factors.

The first four values are:
\begin{eqnarray*}
 r\nm_1 &= &1,\\
 r\nm_2 &= &1/2,\\
 r\nm_3 &= &(q^2+1)^2 (q - 1)/6(q^5-1),\\
 r\nm_4 &= &h\nm(q-1)^2/24(q^5-1)(q^9-1), 
\end{eqnarray*}
where 
$h\nm = q^{12} - q^{11} + 4q^{10} + 3q^8 + 4q^7 - q^6 + 4q^5 + 3q^4 +4q^2 -q +1.$
By convention, we set $r\nm_0=1$.
We will need the following generalization of Lemma~\ref{lem:2}.
Again it is an immediate consequence of Lemma~\ref{lem:3}.

Let $g$ be a (not necessarily monic) 
polynomial in $k[x]$, with $\deg(g)\le n$.
Write $P_{n,g}$ for the set of polynomials $f$ in $R[x]$ of degree
at most $n$ whose reduction
$\ov{f}$ is $g$. 
The space $P_{n,g}$ is a measurable subset of $P\nm_n(R)$ and therefore 
inherits a topology and a measure.  We renormalize this measure to give 
$P_{n,g}$ total measure~1.

Recall that if $h$ in $k[x]$ is monic, then we defined $P_h$ 
to be the set of monic polynomials $f$ in $R[x]$ such that $\ov{f}=h$.

\begin{lemma}\label{lem:2nm}
Let $g$ and $h$ in $k[x]$ be relatively prime polynomials of degree
$m$ and $n$, respectively, and assume that $h$ is monic.
Let $m'\ge m$. Then the multiplication map
$$
M: P_{m',g}\times P_h\to P_{m'+n,gh}
$$
is a measure-preserving homeomorphism.
\end{lemma}

We will write $k^{*}$ for the nonzero elements of the residue field $k$.
Let $j=j(f)$ be the multiplicity of $\infty$
as a root of $\ov{f}$, by which we mean that
exactly the first $j$ leading coefficients of $\ov{f}$ are zero.
Evidently, $0\le j\le n+1$.  Then if $j\le n$, by Hensel's lemma
$f$ factors uniquely as $f=f^{\inf}f_{\inf}$, where 
$\deg f_{\inf}\le j$, $\deg f^{\inf}=n-j$, $f_{\inf}$ reduces mod $\pi R$ to 
a nonzero constant, and $f^{\inf}$ is monic.

\begin{lemma}\label{lem:1nm}
If $0\le i\le n$, then:
\begin{enumerate}

\item $\pr(f \spc|j(f)=i)=$
$$
\begin{array}{l}
\pr(f_{\inf}\spc|
\hbox{$\deg f_{\inf}\le i$, $\ov{f_{\inf}}\in k^{*}$}
)\cdot\\
\pr(f^{\inf}\spc|
\hbox{$\deg f^{\inf}=n-i$, $f^{\inf}$\mnc}
).
\end{array}
$$

\item $\pr(f_{\inf}\spc|
\hbox{$\deg f_{\inf}\le i$, $\ov{f_{\inf}}\in k^{*}$}
)=q^{-\binom{i}{2}}r_{i}$.
\item $\pr(f^{\inf}\spc|
\hbox{$\deg f^{\inf}=n-i$, $f^{\inf}$\mnc}
)=r_{n-i}$.
\item
If $i=n+1$ (i.e., if $\ov{f}=0$), then
$$
\pr(f \spc|j(f)=i)=r_n\nm.
$$
\end{enumerate}
\end{lemma}
\begin{proof}
Part~(1) of the lemma follows from Lemma~\ref{lem:2nm}.

After dividing $x^i f_{\inf}(1/x)$ by its leading coefficient, we get a
random monic polynomial reducing to $x^i$.  Applying part~(2) of
Lemma \ref{lem:1} then proves~(2).  

(3) is immediate by definition.

Conditioned on 
$\ov{f}=0$, $f/\pi$ is a random polynomial in $P_n\nm(R)$. This proves 
(4).
\end{proof}

Now write $f=\sum_{0\le i\le n}a_ix^{n-i}$.
Then $j(f)$ is $j$ if and only if
$\ov{a_0}=\cdots=\ov{a_{j-1}}=0$ and (if $j\le n$) $\ov{a_{j}}\ne0$.

Thus, the probability that $j(f)$ equals $j$ is
$\frac{q-1}{q^{j+1}}$ if $0\le j\le n$, and $\frac1{q^{n+1}}$ if $j=n+1$.

Thus, we have
\begin{equation}
 r\nm_n= \sum_{0\le j\le n}\frac{q-1}{q}r_{n-j}\frac{r_j}{q^{\binom{j+1}{2}}}
+ \frac{r\nm_n}{q^{n+1}}.
\end{equation}

This is easily seen to be equivalent to 
$$
\FF\nm = \frac{q-1}{q}\FF\GG,
$$
where 
\begin{eqnarray}\label{fnm}
\FF\nm(t) &= &\sum_{n\ge0}(1-q^{-n-1})r\nm_n t^n, \mbox{ and}\\
\GG(t) &= &\sum_{n\ge0}\frac{r_{n}}{q^{\binom{n +1}{2}}}t^n.\label{g}
\end{eqnarray}

Here is another equivalent formulation:
\begin{corollary}
$\FF\nm = \frac{q-1}{q}\GG^{q+1}.$
\end{corollary}

Using Euler's method, we get
\begin{equation}
\sum_{0\le j\le n}(n-(q+2)j)(1-q^{-(n-j)-1})r\nm_{n-j}
\frac{r_{j}}{q^{\binom{j +1}{2}}}=0.
\end{equation}
In particular, $r_n\nm$ is a rational function of $q$.

Here is a slightly different viewpoint, which leads to a slightly different
formulation of the result.

Let $\PP(k)$ be the set of lines 
through the origin in the plane $k^2$.
Such a line is uniquely determined by its slope, which is 
an element of the set $k$ or undefined.
This identifies $\PP(k)$ with $k\cup\{\infty\}$;
hence, $\PP(k)$ has cardinality $q+1$.

Assume for the moment that $f$ splits completely, and that  $\ov{f}$,
its reduction mod $\pi R$, is nonzero. 
Write
$$
\ov{f} = (\ov{b}_1x-\ov{a}_1)\cdots(\ov{b}_nx-\ov{a}_n).
$$
Then, for each $i$, $(\ov{a}_i,\ov{b}_i)$ in $k^2$ is not the origin, so 
the line through this point and the origin is in $\PP(k)$.
As $i$ varies, we get a a set 
of $n$ points in $\PP(k)$, with multiplicities.

For $\ell\in \PP(k)$, let $d_{\ell}$ be the number of $1\le i\le n$ such that
$(\ov{a}_i,\ov{b}_i)$ lies on $\ell$.
We have $n=\sum_{\ell\in \PP(k)}d_{\ell}$.

For any function $d'$ from $\PP(k)$ to the set $\N$ of nonnegative integers,
set  $\abs{d'}=\sum_{\ell\in \PP(k)}d'_{\ell}$.
Reasoning as in Lemmas \ref{lem:1} and \ref{lem:1nm} now yields
\begin{theorem}\label{nm}
Let $R$ be a complete discrete valuation ring with finite residue
field $k$.  Set $q=\abs{k}$. For $n$ a positive integer, write $r\nm_n$ for
the probability that a random degree $\le n$ polynomial over $R$
splits into linear factors. Set $r\nm_0=1$. Then 
$$
r\nm_n=\frac{r\nm_n}{q^{n+1}} + 
\frac{q-1}{q}\sum_{\abs{d}=n}\ \prod_{i\in\PP(k)}
\frac{r_{d_i}}{q^{\binom{d_i +1}{2}}}.
$$
\end{theorem}

\section{Generating functions}
In this section we regard $q$ as a variable, and define the sequence
of rational functions $(r_n)=(r_n(q))$ by $r_0=r_1=1$,
and by \begin{eqnarray}
\FF(t) &=&\sum_{n\ge0}r_n t^n, \label{def1} \\
\GG(t) &=&\sum_{n\ge0}\frac{r_n}{q^{\binom{n +1}{2}}}t^n, \qquad \hbox{and} \label{def2}\\
\FF(t) &=& \exp(q \log \GG(t))=\GG(t)^q.\label{def3}
\end{eqnarray}

Thus, if we plug in a prime power for $q$ we recover the $r_n$'s with their
previous meaning.

We note some properties of these rational functions.

\begin{lemma}\label{powerseries}\ 

\begin{enumerate}
\item The degree of the numerator of $r_n$ is the degree of the denominator,
and $\lim_{q\rightarrow\infty} r_n(q) = 1/n!$\ .
\item  $r_n$ vanishes at 0 to order $\binom n2$.
\item The only poles of $r_n$  are at roots of unity.
\item $r_n(q) = r_n(1/q)  q^{\binom n 2}$.
\item If $q$ is fixed, $|q|\ge 3$ and $q$ is not an integer then,
as $n\rightarrow\infty$,
$$
r_n(q)\sim (-z_0)^q c^q \Gamma(-q)^{-1} n^{-q-1} z_0^{-n}
$$
for some complex numbers $z_0$, $c\ne0$.
\end{enumerate}
\end{lemma}

{\it Proof of Properties {\rm(1)--(4)}.}
Property (1) follows from Observation~3 in Section~\ref{secff}.

For~(4), take the defining equations (\ref{def1}), (\ref{def2}) and
(\ref{def3}) of $r_n$, and substitute $qt$ for $t$.  We get
$$
\sum_{n\ge0} r_n(q) q^n t^n=\left(\sum_{n \ge0}
\frac{r_n(q)}{q^{\binom{n+1}{2}}} q^n t^n\right)^q,
$$
which we can rewrite as
$$
\sum_{n\ge0} \frac{r_n(q) q^n}{q^{\binom{n+1}{2}}}
\frac{1}{(1/q)^{\binom{n+1}{2}}} t^n=
 \left(\sum_{n \ge0}
\frac{r_n(q)q^n}{q^{\binom{n+1}{2}}} t^n\right)^q,
$$
or as
$$
\left(\sum_{n\ge0} \frac{r_n(q)}{q^{\binom{n}{2}}}
\frac{1}{(1/q)^{\binom{n+1}{2}}} t^n\right)^{1/q}=
 \sum_{n \ge0}
\frac{r_n(q)}{q^{\binom{n}{2}}} t^n.
$$
But this is the defining equation for $r_n(1/q)$, with $r_n(q)/q^{\binom{n}{2}}$
substituted for $r_n(1/q)$.  Also, for $n=0$ and $n=1$, 
$r_n(q)/q^{\binom{n}{2}}=1$.  Therefore we must have 
$r_n(1/q)=r_n(q)/q^{\binom n 2}$.

For (3), 
define $\bar s_n=\bar s_n(q)$ by $\bar s_n(q)=r_n(q)q^{-\binom n2}$.
We will prove by induction on $n$ that the only poles of $\bar s_n$ are at 
roots of unity.  Since $\bar s_0=\bar s_1=1$, this is clear for $n<2$.  If 
$n\ge 2$, multiply Corollary \ref{cor:1} by $q^n$ and rewrite it as
$$
\sum_{0\le j\le n} (n-(q+1)j) \bar s_{n-j} \bar s_j q^{\binom{n-j+1}{2}}=0.
$$
Separating out the terms involving $\bar s_n$ gives
$$
n q^{\binom{n+1}{2}} \bar s_n - nq \bar s_n=qt_n,
$$
where $t_n$ is a polynomial in $q$, $\bar s_1$, \dots, $\bar s_{n-1}$.
But then $\bar s_n=t_n/n(q^{\binom{n+1}{2}-1}-1)$, so the result follows from
the induction hypothesis.

Property~(2) is an immediate consequence of (1) and (4).
This finishes the proof of Properties {\rm(1)--(4)}.\qed
\medskip

Property~(5) will be proved in the remainder of the section.
By property~(3), it makes sense to plug in for $q$ any complex number
in the complement of the unit disk, $\mathcal{S}=\{q|\,\abs{q}>1\}$. 
For $q\in{\mathcal{S}}$ 
not a positive integer, we expect that the power series $\GG$ 
has a complex zero $z_0$ (necessarily nonzero, as $r_0=1$)
which is simple and is the 
unique zero in the disc $\{z|\,|z|\le |z_0|\}$).
In fact, this is the case if $|q|\ge 3$ by Theorem~\ref{thm:m} below.
Property (5) follows from this theorem by singularity analysis, with
the constant $c=\GG'(z_0).$

Our particular focus is on the asymptotics of $r_n(q)$, as $n\to\infty$.
Surprisingly, the behavior in $\cal{S} \cap \N$ is 
radically different from the behavior in $\cal{S} \setminus \N$.

To explain this, for $q\in\cal{S}$, we 
note that the terms do not grow too fast.

\begin{lemma}\label{powseries2}
For each $q\in\cal{S}$, there is some positive real $M=M(q)$ such that
$|r_n(q)|\le M^n$ for all $n$.  
\end{lemma}
\begin{proof}
Since $|q|>1$, the series $\sum_{n\ge 1} |q|^{1-\binom{n+1}{2}}$ converges,
so pick $k$ such that $\sum_{n\ge k} |q|^{1-\binom{n+1}{2}}< 1$, 
and then pick $0<\epsilon<1$ such that 
\begin{equation}
\label{eyuni}
\epsilon\sum_{1\le n<k} |q|^{1-\binom{n+1}{2}}<
1-\sum_{n\ge k} |q|^{1-\binom{n+1}{2}}.
\end{equation}
Finally, pick $M\ge 1$ large enough so that $|r_n|<\epsilon M^n$ for $n=1$, 
\dots, $k-1$.  We will prove by induction that for all $n$, $|r_n|\le M^n$.
We may assume that $n\ge \max(k,2)$.  Now for $n\ge 2$, Corollary 
\ref{cor:1} tells us that
$$
r_n=
\frac{\sum_{1\le j\le n-1} (n-(q+1)j) r_{n-j} r_jq^{-\binom{j+1}{2}}
}
{n(q^{1-\binom{n+1}{2}}-1)}.
$$
However, $|n-(q+1)j|\le n-j+|q|j = n+(|q|-1)j\le n|q|$, so the 
absolute value of the numerator is bounded by
$$
\sum_{1\le j\le n-1} n|r_{n-j}||r_j||q|^{1-\binom{j+1}{2}} 
$$
and by our assumption and the induction hypothesis this is less than
$$
n M^n 
(\epsilon \sum_{1\le j< k} |q|^{1-\binom{j+1}{2}}
+ \sum_{k\le j\le n-1} |q|^{1-\binom{j+1}{2}}).
$$
Now the absolute value of the denominator is no smaller than
$n(1-|q|^{1-\binom{n+1}{2}})$, so to have $|r_n|\le M^n$ it will do to have
$$
\epsilon \sum_{1\le j< k} |q|^{1-\binom{j+1}{2}}
+ \sum_{k\le j\le n-1} |q|^{1-\binom{j+1}{2}}
\le 1-|q|^{1-\binom{n+1}{2}},
$$
which follows immediately from (\ref{eyuni}).  This completes the proof.
\end{proof}

From this it follows that 
the series $\GG$ has an infinite radius of convergence.
If $q>1$ is an integer, then $\FF=\GG^q$ implies that $\FF$ 
also has an infinite
radius of convergence, whence its coefficients are small:
$r_n(q)=O(\epsilon^n)$ for all $\epsilon>0$.
(In fact, much more is true, as we shall see in the next section.)
For $q$ nonintegral, on the other hand, property (5) says that this will not
in general be true.

\begin{lemma}\label{powseries3}
If $|q|\ge 2$, $|r_n(q)|\le 1$ for all $n$.
\end{lemma}
\begin{proof}
We proceed by induction on $n$.  
For $n=0$ and $n=1$ the result is trivial. So let 
$n\ge 2$.  Corollary \ref{cor:1} then gives the bound
$$
|r_n|\le
\frac{
\sum_{1\le j\le n-1} (n-j+|q|j) |r_{n-j}| |r_j| |q|^{-\binom{j+1}{2}}
}
{n(1-|q|^{1-\binom{n+1}{2}})},
$$
and by the induction hypothesis, this is no bigger than
$$
\frac{
\sum_{1\le j\le n-1} (n-j+|q|j) |q|^{-\binom{j+1}{2}}
}
{n(1-|q|^{1-\binom{n+1}{2}})},
$$
so it will do to show that
$$
\sum_{1\le j\le n} (n-j+|q|j) |q|^{-\binom{j+1}{2}} \le n.
$$
Since the left-hand side of this inequality is nonincreasing in $|q|$,
it will do to show that
\begin{equation*}
\sum_{1\le j\le n} (n+j) 2^{-\binom{j+1}{2}} \le n.
\end{equation*}
For $n=2$ this can be directly verified.  Otherwise, $n\ge 3$.  Now observe 
that 
$$
\sum_{1\le j\le n} n 2^{-\binom{j+1}{2}}\le 
\sum_{j\ge 1} n 2^{-\binom{j+1}{2}}
\le \sum_{j\ge 1} n 2^{-(2j-1)} = \frac{2n}{3}
$$
and
$$
\sum_{1\le j\le n} j 2^{-\binom{j+1}{2}}
\le \sum_{j\ge 1} j 2^{-(2j-1)} = \frac{8}{9}\le 1\le \frac{n}{3}.
$$
This completes the proof.
\end{proof}

\begin{theorem}\label{thm:m}
If $q=2$ or $|q|\ge 3$, 
$\GG(z)$ has a unique zero in the disc $D=\{z|\,|z|<|q|+1\}$, which is simple.
\end{theorem}
\begin{proof}
Set $\GG_0(z)=r_0+r_1z/q+r_2z^2/q^3=1+z/q+z^2/(2q^2(q+1))$.  
Then $\GG_0$ has roots at 
$-(q+q^2)\pm q^2\sqrt{1-q^{-2}}$.  By looking at the power series for 
$\sqrt{1+x}$, we can write
$\sqrt{1-q^{-2}}=1-q^{-2}/2+\epsilon$, where $|\epsilon|<|q|^{-4}/5$.  
Plugging in
this estimate yields roots $r_1=-q-1/2+q^2\epsilon$ and 
$r_2=-2q^2-q+1/2-q^2\epsilon$, and since 
$$
|r_1|\le|q|+\frac12+
\frac{|q|^{-2}}{5}
\le |q|+\frac{11}{20}<|q|+1
$$
and 
$$
|r_2|>2|q|^2-|q|-\frac12-\frac{|q|^{-2}}{5}
\ge 2|q|^2-|q|-\frac{11}{20}>|q|+1,
$$ 
exactly one of these roots
is in $D$.  By Rouch\'e's Theorem, it will do to show that 
$|\GG(z)-\GG_0(z)|<|\GG_0(z)|$ on $C$, 
where $C=\{z|\,|z|=|q|+1\}$.  Let $z\in C$.
We then have
\begin{eqnarray}
|(q+1)\GG_0(z)|&\ge&\frac12|q^{-2}|(|z|-|r_1|)(|r_2|-|z|)\label{eralph0}\\
&\ge&\frac12|q|^{-2}\frac{9}{20}(2|q|^2-2|q|-\frac{31}{20})\nonumber\\
&\ge&\frac{9}{20}(1-|q|^{-1}-\frac{31|q|^{-2}}{40})\nonumber\\
&\ge&\frac{441}{3200}.\nonumber
\end{eqnarray}

Now,
\begin{equation}
\label{eralph}
|(q+1)r_3\frac{z^3}{q^6}|=
\frac16(1+\frac{1}{|q|})^3\frac{|q^2-q+1||q-1|}{|q^5-1|}.
\end{equation}
If $q=2$, the right-hand side of (\ref{eralph}) is $27/496$,
and if $|q|\ge 3$, it is no more than
$$
\frac16(1+\frac{1}{|q|})^3|q|^{-2}
\frac{(1+|q|^{-1}+|q|^{-2})(1+|q|^{-1})}{1-|q|^{-5}}
\le
\frac16(\frac43)^3 3^{-2}\frac{(13/9)(4/3)}{1-3^{-5}}=
\frac{832}{9801}.
$$
In any case, it is no more than $9/100$.
Again, we have
\begin{equation}
\label{eralph2}
|(q+1)r_4\frac{z^4}{q^{10}}|=
\frac1{24}(1+\frac{1}{|q|})^4\frac{|h(q)||q-1|^2}{|q+1||q^5-1||q^9-1|}.
\end{equation}
If $q=2$, the right-hand side of (\ref{eralph2}) is 
$1161/2027648$.
If $|q|\ge 3$, $|q^6h(q)|\le 3|(q^5-1)(q^9-1)|$, so it is no more than
$$
\frac1{24}(1+\frac{1}{|q|})^4|q|^{-5}
\frac{3(1+|q|^{-1})^2}{1-|q|^{-1}}
\le
\frac1{24}(\frac{4}{3})^43^{-5}\frac{3(1+3^{-1})^2}{1-3^{-1}}
=\frac{256}{59049}.
$$
In any case, it is no more than $1/200$.
Finally, using Lemma \ref{powseries3}, we have
\begin{eqnarray*}
\left|
\sum_{n\ge 5}
\frac{(q+1)r_n}{q^{n(n+1)/2}}z^n\right|
&\le&|1+q^{-1}|\sum_{n\ge 5} |\frac zq|^n |q|^{-(n(n-1)/2-1)}\\
&\le&(1+|q|^{-1})\sum_{n\ge 5} (1+|q|^{-1})^n |q|^{-(n(n-1)/2-1)}\\
&\le& \frac32\sum_{n\ge 5} (\frac32)^n 2^{-(n(n-1)/2-1)}\\
&\le& (\frac32)^6 2^{-9} \frac1{1-\frac32 2^{-5}}\\
&=& \frac{729}{31232}\le \frac{1}{40}.
\end{eqnarray*}
Adding up these inequalities tells us that
$$
|(q+1)(\GG_0(z)-\GG(z))|=
\left|
\sum_{n\ge 3}
\frac{(q+1)r_n}{q^{n(n+1)/2}}z^n\right|
\le \frac{9}{100}+\frac{1}{200}+\frac{1}{40}=\frac{3}{25}.
$$
Together with (\ref{eralph0}), this proves the result.
\end{proof}

As remarked above, property (5) follows from Theorem~\ref{thm:m}.

\section{Open questions}
We mention some open questions for further research.

{\bf 1.  Zeroes of the numerator.}\ Let $N_n$ be the numerator of $r_n$.
Is it true that $N_n$ has no zeroes on the negative real axis?

{\bf 2.  Zeroes of $\GG$.}\ Let $q\ge2$ be an integer. 
Is it true that the zeroes of the function $\GG$ are 
all real and negative?

We remark that the zeroes $|z_0|\le|z_1|\le\cdots$ of $\GG$
completely determine $\GG$ as follows.
Because $\GG(t)$ is an entire function of order zero, i.e.,
for all $\alpha>0$, $|\GG(z)|=o(e^{|z|^\alpha})$ as $|z|\to\infty$,
it is equal up to a constant factor to the Weierstrass product 
$\prod_{i\ge0} (1 - \frac{t}{z_i})$. But then $\GG(t)$ equals this
product since $\GG(0)=1$.

{\bf 3. Asymptotics.}\ Prove the more accurate asymptotic expansion
(\ref{zzz}) for $\log_q r_n(q)$, for $q>2$ an integer.
\section{Asymptotics}
\label{secasymp}
In this section $q\ge2$ is an integer.
If $q=2$, we give the first three terms in the asymptotic expansion of 
$\log_2 r_n$.  If $q>2$, we give the first two terms in the 
asymptotic expansion of $\log_q r_n$ and a guess for its third term.

\begin{theorem}\label{thm:moews}
If $q=2$, then there is a constant $K_2$ and a 
continuous function $\bar W_2$ with period 1 such that
\begin{equation}
\label{e77}
\abs{\log_2 r_n + \frac{n^2}{2}  + \frac12 n\log_2 n
-(\frac{1}{2}+\bar W_2(\log_2 n))n}\le K_2,
\end{equation}
for all $n\ge 1$.
(Figure~\ref{fig2} shows a graph of $\bar W_2$.)
If $q\ge2$ is an integer, then there is a constant $C_q$ such that 
\begin{equation}
\label{e78}
\abs{\log_q r_n + \frac{n^2}{2(q-1)}  + \frac12 n\log_q n}\le C_q n,
\end{equation}
for all $n\ge1$, and we conjecture (Conjecture \ref{c26}) that
there is a constant $K_q$ and a continuous function $\bar W_q$ with
period 1 such that
$$
\abs{\log_q r_n + \frac{n^2}{2(q-1)}  + \frac12 n\log_q n
-(\frac{1}{2}+\bar W_q(\log_q n))n}\le K_q,
$$
for all $n\ge 1$.
\end{theorem}

The idea of the proof of \eqref{e78} is this. In the proof of the 
exact formula for $r_n$, we looked at a certain infinite tree
whose branching was given by the factorization of a polynomial 
mod the maximal ideal $\pi R$ of $R$, mod $\pi^2R$, etc.
In a modification of this method,  
we use finite trees instead.
Lemma~\ref{lem:trees} 
expresses $r_n/q^{\binom{n+1}{2}}$ as the sum of a certain function $H$
over labelled $q$-trees (defined below) with $n$ leaves.
By Lemma~\ref{blah}, we will see that for any fixed $q$, 
the logarithm of the number of labelled $q$-trees is $O(n)$.
For example, consider the case $q=2$.  In this case our labelled $q$-trees 
are the same as binary trees.
The number of binary trees with
$n$ leaves is given by the Catalan number $\frac{1}{n}\binom{2n-2}{n-1}.$
By Stirling's formula, this has logarithm $O(n).$

If we set $M$ to be the maximum value of $\log_q H$, we then have
$$
M\le-{\binom{n+1}{2}}+ \log_q r_n\le M+O(n).
$$

Thus, to prove \eqref{e78},
it will be enough to estimate
the maximum value of $\log_q H$ to within $O(n)$.
We will see by Lemma~\ref{bal} that 
$H$ is maximized at a tree we will call the {\em well-balanced $q$-tree},
and in Subsection~\ref{ss:moews} we will calculate $H$ for the
well-balanced $q$-tree.

Although \eqref{e77} is a refinement of \eqref{e78}, its proof is
different, and relies on a direct use of \eqref{maineq}, together with an
estimate of how much the largest term appearing in \eqref{maineq}
contributes to \eqref{maineq}.

\subsection{$q$-trees}
A rooted tree is a connected acyclic graph with a distinguished vertex
(the root).
We can direct the edges of a rooted tree in a unique way,
by directing them away from the root.
The root has in-degree 0; all other vertices have in-degree~1.
The edges emanating from a vertex 
go to distinct vertices, called the {\em children}
of $v$. A vertex with no children is a \emph{leaf} vertex.

A \emph{$q$-tree} is a rooted tree
in which the number of out-edges from each vertex 
is at most $q$, and is not 1 (i.e. there is branching at
each vertex that is not a leaf.)

For $v$ a vertex of the $q$-tree $T$, we write $T_v$ for the
full subtree whose vertices are
$v$ and all its descendants. This is a $q$-tree.

A \emph{labelled $q$-tree} is a $q$-tree together with a labelling of its
edges with elements of the set $S =\{0,1,\dots,q-1\}$ so that 
the out-edges emanating from each 
vertex have distinct labels.

\begin{example}We list the $q$-trees with at most $3$ leaves. 
\begin{enumerate}
\item[(i)] The tree $\tau_1$ with one vertex and no edges is the unique
$q$-tree with one leaf. In all other $q$-trees, a vertex 
is a leaf if and only if it has total degree~1.
\item[(ii)] The tree $\tau_2$ with three vertices, a root with two
leaf children, is the unique $q$-tree with two leaves.
There are $\binom{q}{2}$ ways of labelling this $q$-tree.
\item[(iii)] There are two $q$-trees with three leaves (one if $q=2$). 
One, call it $\tau_{3a}$,
has a root with three leaf children (if $q>2$). The other,
$\tau_{3b}$, has
a root with two children $v$ and $w$, with $v$ a leaf and $T_w=\tau_2$.
There are $\binom{q}{3}$ ways of labelling $\tau_{3a}$, and 
$q(q-1)\binom{q}{2}$ ways of labelling $\tau_{3b}$.
\end{enumerate}
\end{example}

\begin{lemma}\label{blah}
The number of labelled $q$-trees with $l$ leaves is at most 
$(2q+1)^{5l-3}$.
\end{lemma}

\begin{proof}
Let $\Sigma$ be the set of $2q+1$ symbols: the parentheses 
$(_i$ and $)_i$ for $0\le i\le q-1$, 
plus a dot.
We construct for every labelled $q$-tree a distinct string
of at most $5l-4$ of these symbols, thus constructing
an injective map from the set of labelled $q$-trees with $l$ leaves to 
$\cup_{0\le i\le 5l-4} \Sigma^i.$

We proceed by induction on $l$.
For $l=1$, the tree $\tau_1$ corresponds to the dot.

Assume $l\ge2$, and 
let $T$ be a labelled $q$-tree with $l$ leaves.
Let the root of $T$ have $j$ children.
For each child $v$ of the root of $T$,
bracket the sequence corresponding to the 
subtree $T_v$ with the symbols $(_i$ and $)_i$, where 
$i\in S$ is the label of the edge from the root to the vertex $v$.
By the induction hypothesis, this uses at most
$5l-4j+2j=5l-2j\le 5l-4$ symbols.

The number of labelled $q$-trees with $l$ leaves
is therefore at most 
$$
\sum_{0\le i\le 5l-4}(2q+1)^i
=(2q)^{-1}((2q+1)^{5l-3}-1)\le(2q+1)^{5l-3}.
$$
\end{proof}

\subsection{A $q$-tree recursion}\label{subsec:rec}
Equating coefficients of $t^n$ in the defining equations
(\ref{def1}), (\ref{def2}) and (\ref{def3}) of $r_n$ gives us
\begin{equation}
r_n = \sum_{\abs{b}=n}\prod_{i=0}^{q-1}\frac{r_{b_i}}{q^{\binom{b_i+1}{2}}}.
\end{equation}

Group the $r_n$ terms to get

\begin{equation*}
r_n = \frac{q}{q^{\binom{n+1}{2}}}r_n + 
{\sum_{\abs{b}=n}}'\prod_{i=0}^{q-1}\frac{r_{b_i}}{q^{\binom{b_i+1}{2}}},
\end{equation*}
where $\sum'$ is the sum over $b$ such that $b_i>0$ for at least two values
of $i$.

Setting $s_n=r_n/{q^{\binom{n+1}{2}}}$, we have

\begin{equation*}\label{m-}
q^{\binom{n+1}{2}}s_n = qs_n + 
{\sum_{\abs{b}=n}}'\prod_{i=0}^{q-1}s_{b_i}.
\end{equation*}

For $n\ge2$, set $\bt_n=1/{(q^{\binom{n+1}{2}}-q)}$.
We have proved the following.

\begin{lemma}\label{m}
Assume $n\ge2$. Then, 
\begin{equation}
\label{eqm}
s_n = \bt_n{\sum_{\abs{b}=n}}'\prod_{i=0}^{q-1}s_{b_i}.
\end{equation}
\end{lemma}

\begin{remark}\label{remark:t}
We can decompose a $q$-tree into its root, plus a subtree $T_v$ for each
child $v$ of the root of $T$. This decomposition underlies
Lemma~\ref{blah}, and will allow us to interpret
Lemma~\ref{m} as a recursion on labelled $q$-trees.
\end{remark}

Write $\ell(v)$ for the number of leaves of the tree $T_v$.
As an easy application of the decomposition in Remark~\ref{remark:t},
we have the following. Let $T$ be a $q$-tree with more than one vertex.
Then the root of $T$ is not a leaf, and the number of leaves of $T$ equals
$\sum\ell(v)$, where the sum is over the set of
children $v$ of the root of the tree $T$.

We can now interpret the recursion for $s_n$ in terms of labelled $q$-trees.
Set $\bt_1=s_1=1/q$.  
\begin{lemma}\label{lem:trees}We have
\begin{equation}
\label{eqn:trees}
s_n=\sum \prod_{v\in T}\bt_{\ell(v)},
\end{equation}
where the sum is over labelled $q$-trees $T$ with $n$ leaves, and 
the product is over all vertices $v$ of the tree $T$.
\end{lemma}
\begin{proof}
For $n=1$, this follows by our choice of $\bt_1$.
The general case follows by Remark~\ref{remark:t} and Lemma~\ref{m}.
\end{proof}
We shall write $H(T)$ for $\prod_{v\in T}\bt_{\ell(v)}.$ 

\begin{example}\label{ex:trees}
We recalculate $r_1,r_2$ and $r_3$ using Lemma \ref{lem:trees}.
\begin{enumerate}
\item[(i)] We have $s_1=1/q$ , so that $r_1=1$.
\item[(ii)] We have $H(\tau_2)=\bt_1^2\bt_2$.
There are $\binom{q}{2}$ ways of labelling the edges of $\tau_2$.
Hence $r_2=q^3s_2=q^3\binom{q}{2}\bt_1^2\bt_2=q/{2(q+1)}$.
\item[(iii)] We have $H(\tau_{3a})=\bt_3\bt_1^3$.
There are $\binom{q}{3}$ ways of labelling the edges of 
$\tau_{3a}$ (even if $q=2$!). 
The contribution to $r_3$ from $\tau_{3a}$ is thus 
$q^6\binom{q}{3}\bt_3\bt_1^3$. 
We have $H(\tau_{3b})=\bt_3\bt_2\bt_1^3$.
There are $q(q-1)\binom{q}{2}$ ways of labelling the edges of 
$\tau_{3b}$. Hence, the contribution to $r_3$ from $\tau_{3b}$ is  
$q^6q(q-1)\binom{q}{2}\bt_3\bt_2\bt_1^3$. Add, and this agrees with
the formula in Section~\ref{intro}.
\end{enumerate}
\end{example}

\subsection{The well-balanced $q$-tree}
In Subsection~\ref{subsec:rec} we defined $\bt_n=1/(q^{\binom{n+1}2}-q)$ 
for $n\ge2$ and $\bt_1=1/q$.

\begin{lemma}\label{beta}We have 
\begin{description}
\item[\rm (i)] $\log_q \bt_n =-\binom{n+1}{2} + o(1)$ as $n\to\infty$.
\item[\rm (ii)] The sequence $(\bt_n/\bt_{n-1})_{n=2}^{\infty}$ 
is monotone decreasing.
\end{description}
\end{lemma}
\begin{proof}
Part (i) is clear since $\bt_nq^{\binom{n+1}{2}}$ tends to~1
as $n\to\infty$.

In order to prove (ii) it suffices to show 
$\frac{1}{\bt_n^2}\le \frac1{\bt_{n-1}\bt_{n+1}}$ for all $n\ge2$.

First consider the case $n=2$. We have
$\frac{1}{\bt_2^2}=(q^3-q)^2\le(q^3-q)(q^3+q)= q^6- q^2\le q(q^6- q),$
and this last quantity is equal to $1/\bt_{1}\bt_{3}.$

Next suppose $n\ge3$. Then $\frac{1}{\bt_n^2}=(q^{\binom{n+1}{2}}-q)^2$,
and we have  
\begin{eqnarray*}
(q^{\binom{n+1}{2}}-q)^2&\le &(q^{\binom{n+1}{2}}-q)
(q^{\binom{n+1}{2}}+q)\\
&=&q^{n^2+n}-q^2\\
&\le &(q^{\binom{n}{2}-1}+q)(q^{\binom{n+2}{2}}-q)=AB, \qquad \hbox{say}.
\end{eqnarray*}
The last inequality follows since $\binom{n}{2}-1$ is less than
$\binom{n+2}{2}$ and their sum is $n^2+n$.

Since $2\le q,$ we have $2\le q(q-1)$, which implies that
$2q\le q^j(q-1)$ for any $j\ge2$. Rearranging the terms gives
$q^{j}+q\le q^{j+1}-q,$ which implies (taking 
$j=\binom{n}{2}-1\ge\binom32-1=2$) 
that $A\le 1/\bt_{n-1}$. As $B=1/\bt_{n+1}$, the lemma is proved.
\end{proof}

Let $n\ge 1$ be an integer.
We are interested in $q$-tuples of nonnegative integers
that sum to $n$, and such that each pair differs by at most one.
We remark that such a $q$-tuple $i$ exists: write $n=qx+y$ with 
$0\le y\le q$ (we allow either $y=0$ or $y=q$); now take 
$i_1=\cdots=i_y=x+1$ and $i_{y+1}=\cdots=i_q=x$.
Moreover, if such a $q$-tuple contained a value less than $x$
(resp.\ larger than $x+1$), 
all values would be at most $x$ (resp.\ at least $x$),
and the sum would be too small
(resp.\ too large). Thus all values are $x$ or $x+1$, and the number of each 
must be $q-y$ (resp.\ $y$).
Hence the $q$-tuple is unique up to order. 

\begin{lemma}
\label{ilovemevoli}
Let $n\ge 1$ be an integer.
There is a unique $q$-tree $T=T(n)$ with $n$ leaves such that 
for every vertex $v$ of $T$:
\begin{enumerate}
\item If $\ell(v)<q$ then all children of $v$ are leaves.
\item If $\ell(v)\ge q$ then $v$ has $q$ children and,
for any two children $w$ and $w'$ of $v$,
$\ell(w)$ and $\ell(w')$ differ by at most~1.
\end{enumerate}
\end{lemma}
\begin{proof} If $n=1$, then $T$ is the unique $q$-tree with one
(leaf) vertex.  If $1<n<q$, then $T$ is the unique $q$-tree with 
$n+1$ vertices that consists of 
a root with $n$ children all of which are leaves.

For $n\ge q$, we define $T=T(n)$ by induction on $n$.
Write $n=qx+y$ with $0\le y\le q-1$ 
as in the remarks preceding this lemma.
Applying property 2 above at the root of $T$, we see that
if $T$ exists, its root must have children $v_1$, \dots, $v_q$ 
such that $\ell(v_1)=\cdots=\ell(v_y)=x+1$ and $\ell(v_{y+1})=\cdots
=\ell(v_q)=x$.  It now follows from the induction hypothesis that
each $T_{v_i}$ must equal $T(\ell(v_i))$.  This gives a unique candidate
for $T$, namely
the $q$-tree whose root has precisely $q$ children, and such that 
$T_w=T(x+1)$ for $y$ of these children $w$ 
and $T_w=T(x)$ for the remaining $q-y$ children $w$.
It is easy to see that this does indeed satisfy properties 1 and 2.
\end{proof}

We call a $q$-tree {\em well-balanced} if it satisfies the
two conditions of the lemma.  For $T$ the well balanced $q$-tree with
$n\ge 1$ leaves, write $\TTT_n=H(T)=\prod_{v\in T} \bt_{\ell(v)}$.
Set $\TTT_0=1$.

\begin{lemma} \label{balpre} \ 
\begin{enumerate}
\item[\rm (i)] $\TTT_1=1/q$. 
\item[\rm (ii)] If $n\ge 2$ and we write $n=qx+y$, $0\le y\le q$, then
$$
\TTT_n=\bt_n \TTT_x^{q-y} \TTT_{x+1}^y.
$$
\item[\rm (iii)] The sequence $(\TTT_i/\TTT_{i-1})_{i=1}^{\infty}$ 
is monotone decreasing.
\end{enumerate}
\end{lemma}
\begin{proof}  (i) is obvious, and (ii) follows immediately
from the defining properties of the well-balanced $q$-tree and the remarks
preceding Lemma \ref{ilovemevoli}.
We now prove that for all $n$, the sequence
$(\TTT_i/\TTT_{i-1})_{i=1}^n$ is monotone decreasing; this will prove (iii).
We induce on $n$.
If $n\le 3$, the result follows from Example \ref{ex:trees}, so 
assume that $n\ge 4$.
We need to show that $\TTT_{n-1}/\TTT_{n-2}\ge\TTT_{n}/\TTT_{n-1}$.

For any $2\le m\le n-1$, write $m=qx+y$ with $0\le y\le q-1$
(note that we do not allow $y=q$). 
By (ii), $\TTT_m=\bt_m\TTT_x^{q-y}\TTT_{x+1}^y$.

Next we calculate $\TTT_{m+1}$.
We have $m+1=qx+(y+1)$, with $y+1\le q$, so (ii) gives
$\TTT_{m+1}= \bt_{m+1}\TTT_x^{q-y-1}\TTT_{x+1}^{y+1},$ whence
$\TTT_{m+1}/\TTT_m=\frac{\bt_{m+1}}{\bt_m}\frac{\TTT_{x+1}}{\TTT_x}.$

We apply this to $m=n-2$ and $m=n-1$.
The desired result now follows from part (ii) of Lemma~\ref{beta},
after possibly using the induction hypothesis.
\end{proof}

Let $\TT_n$ be the largest tree contribution $H(T)=\prod_{v\in T}\bt_{\ell(v)}$
among $q$-trees $T$ with $n$ leaves.
Note that here the labelling of the trees is irrelevant.
We set $\TT_0=1$.
The following lemma says that the well-balanced $q$-tree always
yields a largest tree contribution.

\begin{lemma}\label{bal} For all $n\ge 0$, $\TT_n=\TTT_n$.
\end{lemma}
\begin{proof}
We proceed by induction on $n$.  For $n=0$ and $n=1$, the result is obvious,
so let $n>1$.  It will do to prove that $\TT_n\le \TTT_n$.

The contribution from a tree with $n$ leaves is at most
$\bt_n\TT_{i_1}\cdots\TT_{i_q}$,
where $i_1$, \dots, $i_q$
are nonnegative integers that sum to $n$, at least two
of which are positive.
Since at least two $i_z$'s are positive, we have $i_1$,\dots, $i_q\le n-1$,
so by the induction hypothesis,
$\TT_{i_z}=\TTT_{i_z}$ for all $z$.  It therefore follows that
the contribution is no more than
\begin{equation}
\label{foof}
\bt_n\TTT_{i_1}\cdots\TTT_{i_q}.
\end{equation} 
Now by part (iii) of Lemma \ref{balpre},
\begin{equation}
\label{neq}
\TTT_{i-1}\TTT_j\le\TTT_{i}\TTT_{j-1}\mbox{ if } 1\le i\le j.
\end{equation}
Let $n=qx+y$ with $0\le y\le q$.
It follows from \eqref{neq} that 
\eqref{foof} is
largest when $\abs{i_z-i_{z'}}\le 1$ for all $1\le z,z'\le q$,
when, by the remarks preceding Lemma \ref{ilovemevoli}, it equals
$\bt_n \TTT_x^{q-y} \TTT_{x+1}^y.$  By Lemma \ref{balpre}, this equals
$\TTT_n$, so we are done.
\end{proof}

\begin{lemma}
\label{balcontents}
Let $T$ be the well-balanced $q$-tree with $n$ leaves, let $k\ge 0$,
and let $n=q^kx+y$, $0\le y<q^k$.  Then:
\begin{enumerate}
\item If $n<q^{k-1}$, there are no vertices at distance $k$ from the 
root of $T$.
\item If $q^{k-1}\le n<2q^{k-1}$, 
there are $2(n-q^{k-1})$ vertices at distance
$k$ from the root of $T$, each of which is a leaf.
\item If $2q^{k-1}\le n<q^k$,
there are $n$ vertices at distance $k$ from the root of $T$, each of which
is a leaf.
\item If $n\ge q^k$,
there are $q^k$ vertices $v'$ at distance $k$ from the root 
of $T$, $y$ with $\ell(v')=x+1$ and $q^k-y$ with $\ell(v')=x$.
\end{enumerate}
\end{lemma}
\begin{proof}
We prove the lemma by induction on $k$.  
If $k=0$, these results are clear.
Otherwise: 

To prove the first claim, suppose that $n<q^{k-1}$.  Then by the 
induction hypothesis,  all the vertices at distance $k-1$ from the root
of $T$ are leaves, and therefore there are no vertices at distance $k$
from the root, as desired.

If $q^{k-1}\le n<2q^{k-1}$, then by the induction
hypothesis, there are $q^{k-1}$ vertices $v'$ at distance $k-1$ from the 
root, $2q^{k-1}-n$ of which are leaves and $n-q^{k-1}$ of which
satisfy $\ell(v')=2$, i.e. have 2 children, both leaves; this gives a
total of $2(n-q^{k-1})$ vertices at distance $k$ from the root, all leaves,
which is the second claim.

If $2q^{k-1}\le n<q^k$, then write $n=q^{k-1}x'+y'$,
$2\le x'\le q-1$, $0\le y'<q^{k-1}$.  By the induction hypothesis,
there are $q^{k-1}$ vertices $v'$ at distance $k-1$ from the root,
$y'$ of which satisfy $\ell(v')=x'+1$, 
and $q^{k-1}-y'$ of which satisfy $\ell(v')=x'$.
By the properties of the well-balanced $q$-tree, the vertices $v'$ for
which $\ell(v')=x'$ must have $x'$ children, all leaves.  Similarly,
the vertices $v'$ for which $\ell(v')=x'+1$ must have $x'+1$ children,
all leaves.
This gives $n$ vertices at distance $k$ from
the root, all leaves.  This proves the third claim.

Finally, if $n\ge q^k$, write $y=y_0 q^{k-1}+y_1$, $0\le y_0<q$, 
$0\le y_1<q^{k-1}$.  We then have $n=q^{k-1}(qx+y_0)+y_1$.
By the induction hypothesis, therefore,
there are $q^{k-1}$ vertices $v'$ at distance $k-1$ from the root,
$q^{k-1}-y_1$ of which satisfy $\ell(v')=qx+y_0$.
Since $x\ge 1$, $qx+y_0\ge q$, so each of these vertices will,
by the properties of the well-balanced $q$-tree, have $y_0$ children
$v''$ satisfying $\ell(v'')=x+1$ and $q-y_0$ children $v''$ satisfying
$\ell(v'')=x$.  Similarly, $y_1$ of the vertices at distance $k-1$
from the root will satisfy $\ell(v')=qx+y_0+1$ and will have
$y_0+1$ children $v''$ satisfying $\ell(v'')=x+1$ and $q-y_0-1$ children
$v''$ satisfying $\ell(v'')=x$.  This gives a total
of $y_0(q^{k-1}-y_1)+(y_0+1)y_1=q^{k-1}y_0+y_1=y$ vertices $v''$ at
distance $k$ from the root with
$\ell(v'')=x+1$ and 
$(q-y_0)(q^{k-1}-y_1)+(q-y_0-1)y_1=(q-y_0)q^{k-1}-y_1=q^k-y$
vertices $v''$ at distance $k$ from the root with $\ell(v'')=x$, as desired.
This proves the final claim.
\end{proof}

\subsection{The largest tree contribution}\label{ss:moews}
We wish to estimate the contribution from the well-balanced $q$-tree.
By the previous subsection, the well-balanced $q$-tree gives the largest 
contribution to the sum \eqref{eqn:trees}.

The bound \eqref{e78} in Theorem~\ref{thm:moews} will follow from:
\begin{lemma} $\log_q \TTT_n$, the logarithm of the 
contribution to \eqref{eqn:trees} from the well-balanced $q$-tree, is
\begin{equation}
\label{wb}
-\frac{n^2}{2(1-q^{-1})} 
-\frac{n\log_q n}{2}+O(n).
\end{equation}
\end{lemma}
\begin{proof}
We first treat the contribution $\eta$ to $\TTT_n$ from vertices at distance 
more than $\log_q n$ from the root. 
By Lemma~\ref{balcontents}, 
there are $O(n)$ such vertices, and each contributes a factor of $\bt_1$
to $\TTT_n=\prod_v \bt_{\ell(v)}$.
Thus $\log_q \eta=O(n)$.

Let $k$ be an integer such that $1\le q^k\le n$. 
Write $\omega_k$ for the contribution to $\TTT_n$ 
from the vertices at distance $k$ from the root.

Write $n=q^kx + y,$ with $0\le y\le q^k-1$.
By Lemma~\ref{balcontents}, 
we have $\omega_k=\bt_x^{q^k-y}\bt_{x+1}^y$, so that
\begin{eqnarray*}
\log_q \omega_k&=&(q^k-y)\log_q \bt_x  + y\log_q \bt_{x+1} \\
&=&-(q^k-y)\binom{x+1}{2} -y\binom{x+2}{2} +O(q^k)\\
\end{eqnarray*}
by part (i) of Lemma~\ref{beta}.
So, up to terms of order $q^k$, we have
\begin{eqnarray*}
\log_q \omega_k&=&-\frac{x+1}2[(q^k-y)x+ y(x+2)]\\
&=&-\frac{x+1}2[n+y]\\
&=&-\frac12(n+y + \frac{n-y}{q^k}(n+y))\\
&=&-\frac12(\frac{n^2}{q^k}+n).\\
\end{eqnarray*}

Now sum over $k$ such that $1\le q^k\le n$. Note that $\sum_k q^k$ is 
$O(n)$. 

We get 
\begin{equation}
\log_q\TTT_n =\sum_{1\le q^k\le n}\log_q \omega_k+\log_q \eta
=-\frac{n^2}{2(1-q^{-1})}-\frac{n\log_q n}{2}+O(n).
\end{equation}
The lemma is proved.
\end{proof}
The bound \eqref{e78} in Theorem~\ref{thm:moews} is now proved.
\section{A recursion}\label{recursion}
We consider a recursion of the form 
\begin{eqnarray}
\label{recursion3}
\Omega_n&=& (q-y)\Omega_x+y\Omega_{x+1}+\omega_n,\\
&&\qquad \hbox{where $n\ge 2$, $n=qx+y$, $0\le y\le q$,}\nonumber \\
\Omega_0&=&0.\nonumber
\end{eqnarray}

We now show how to solve this recursion if the sequence 
$(\omega_n)$ does not grow too rapidly with $n$.

\begin{theorem}Suppose $\omega_n=O(n^\kappa)$ for $0\le \kappa<1$.
\label{qwerty}
The solution to (\ref{recursion3}) has the form
$$
\Omega_n=W(\log_q n)n+O(n^\kappa),
$$
where $W$ is a continuous function with period~1.
\end{theorem}
\begin{proof}
Let $n\ge 1$.
It is clear from (\ref{recursion3}) that, if we set $\omega_1=\Omega_1$,
and let $T$ be the well-balanced $q$-tree with $n$ leaves,
then $\Omega_n$ is the sum over all vertices $v$ of $T$ of $\omega_{\ell(v)}$.
Write $\floor{z}$ for the largest integer which is no more than $z$
and $\{z\}$ for $z-\floor{z}$.  Now define $X(m)$ by

$$
\begin{array}{cccl}
X(m)&=&0, \qquad & \hbox{$m<q^{-1}$;}\\
X(m)&=&2\omega_1(m-q^{-1}), \qquad & \hbox{$q^{-1}\le m<2q^{-1}$;}\\
X(m)&=&\omega_1 m, \qquad & \hbox{$2q^{-1}\le m<1$;}\\
X(m)&=&\omega_{1+\floor{m}} \{m\} + \omega_{\floor{m}} (1-\{m\}), 
\qquad & \hbox{$1\le m$.}
\end{array}
$$

Observe that $X$ is continuous and that $X(m)=O(m^\kappa)$.
Now Lemma~\ref{balcontents}
implies that the contribution to $\Omega_n$ from vertices at distance $k$
from the root of $T$ is $q^k X(n/q^k)$, so
$$
\Omega_n= \sum_{k\ge 0} q^k X(\frac{n}{q^k}).
$$
Note that the summand is zero for $k>Z=\floor{\log_q n}+1$.  Now
\begin{eqnarray}
\label{last0}
\frac{\Omega_n}{n}&=&
\sum_{0\le k\le Z} \frac{X(n/q^k)}{n/q^k}\\
&=&
\sum_{-Z\le k\le 0} \frac{X(nq^k)}{nq^k}\\
&=&
\sum_{k\ge -Z} \frac{X(nq^k)}{nq^k}+O(n^{\kappa-1}),
\qquad \hbox{as $X(m)/m=O(m^{\kappa-1})$.}\nonumber
\end{eqnarray}

Set 
\begin{equation}
\label{last1}
W(m)=\sum_{k\ge -1-\floor{m}} \frac{X(q^{m+k})}{q^{m+k}},
\end{equation}
or equivalently, 
\begin{equation}
\label{last2}
W(m)=\zeta(\{m\}), \qquad
\zeta(x)= \sum_{k\ge -1} \frac{X(q^{x+k})}{q^{x+k}}.
\end{equation}

It is now clear from (\ref{last0}) and
(\ref{last1}) that $\Omega_n=W(\log_q n)n+O(n^\kappa)$.
However, from (\ref{last2}), $X(m)/m=O(m^{\kappa-1})$, and the continuity
of $X$, we see that on $[0,1]$, $\zeta$ is a sum of
continuous functions which are uniformly bounded by a convergent series.
Therefore $\zeta$ is continuous on $[0,1]$.  The continuity of $W$
now follows from the continuity of $\zeta$ on $[0,1]$ and the fact that,
since $X(q^{-1})=0$, $\zeta(0)=\zeta(1)$.  This concludes the proof.
\end{proof}

\section{The third term}
We can rewrite part (ii) of Lemma \ref{balpre} as
\begin{eqnarray}
\label{recursion1}
\log_q\TTT_n&=&\log_q \bt_n
+(q-y)\log_q\TTT_x+y\log_q\TTT_{x+1}, \\
&&\qquad
\hbox{where $n\ge 2$, $n=qx+y$, $0\le y\le q$.}\nonumber
\end{eqnarray}

Empirically, we have observed that a similar recursion appears to
hold for $\log_q s_n$.  This is because the main contribution in (\ref{eqm})
comes when all the $b_i$'s differ by at most~1.  
For $q=2$, we can prove this.

\begin{lemma}
\label{lemr}
If $q=2$ and $n\ge 2$, then $s_{n+1}s_{n-1}\le \frac{1}{2}s_n^2$.
\end{lemma}
\begin{proof}

Set
\begin{equation}
\label{e0}
R_m=\frac{s_{m+1}/s_m}{s_m/s_{m-1}}=
\frac{s_{m+1}s_{m-1}}{s_m^2} \qquad (m\ge 2).
\end{equation}
Then
$$
\frac{s_{m+1}/s_m}{s_{m-(i-1)}/s_{m-i}}=
R_m\cdots R_{m-(i-1)} \qquad (m\ge i+1, \ i\ge 0),
$$
so
\begin{equation}
\label{e1}
\frac{s_{m+j} s_{m-i}}{s_m s_{m+j-i}}=
\frac{s_{m+j}/s_m}{s_{m+j-i}/s_{m-i}}=
\prod_{0\le k\le j-1} (R_{m+k}\cdots R_{m+k-(i-1)})
\qquad (m\ge i+1, \ i,\,j\ge 0).
\end{equation}

Setting $j=i$ in \eqref{e1} gives
\begin{equation}
\label{e2}
\frac{s_{m-i} s_{m+i}}{s_m^2}=
R_{m-(i-1)} R_{m-(i-2)}^2 \cdots R_m^i \cdots R_{m+(i-2)}^2 R_{m+(i-1)}
\qquad (m\ge i+1, \ i\ge 0)
\end{equation}
and setting $j=i+1$ gives
\begin{equation}
\label{e3}
\frac{s_{m-i} s_{m+i+1}}{s_m s_{m+1}}=
R_{m-(i-1)} R_{m-(i-2)}^2 \cdots R_m^i R_{m+1}^i \cdots R_{m+i-1}^2 R_{m+i}
\qquad (m\ge i+1, \ i\ge 0).
\end{equation}

Now, recalling that
$$
s_n = \beta_n (s_1 s_{n-1} + s_2 s_{n-2} + \cdots + s_{n-1} s_1)
\qquad (n\ge 2),
$$
we get, for $m\ge 3$,
$$
R_{2m}=
\frac{ \beta_{2m-1} \beta_{2m+1}}{\beta_{2m}^2}
\cdot
$$
$$
\frac{
\displaystyle
(2 s_{m-1} s_m + 2 s_{m-2} s_{m+1} + \sum_{2\le i\le m-2} 2 s_{m-i-1} s_{m+i})
(2 s_m s_{m+1} + 2 s_{m-1} s_{m+2} + \sum_{2\le i\le m-1} 2 s_{m-i} s_{m+i+1})
}
{
\displaystyle 
(s_m^2 + 2s_{m-1} s_{m+1} + \sum_{2\le i\le m-1} 2 s_{m-i} s_{m+i})^2
},
$$
so
\begin{equation}
\label{even1}
R_{2m} = 4 
\frac{ \beta_{2m-1} \beta_{2m+1}}{\beta_{2m}^2}
R_m \cdot 
\end{equation}
$$
\frac{
\begin{array}{c}
\displaystyle
(1 + R_{m-1} R_m + \sum_{2\le i\le m-2} R_{m-i}\cdots R_{m-1}^i R_m^i \cdots
R_{m+(i-1)})\cdot \\
\displaystyle
(1 + R_m R_{m+1} + \sum_{2\le i\le m-1} R_{m-(i-1)}\cdots R_m^i 
R_{m+1}^i \cdots R_{m+i})
\end{array}
}
{
\displaystyle
(
1 + 2 R_m + \sum_{2\le i\le m-1} 2 R_{m-(i-1)} \cdots R_m^i \cdots R_{m+(i-1)}
)^2
}.
$$

Similarly, for $m\ge 2$,
$$
R_{2m+1}=
\frac{ \beta_{2m} \beta_{2m+2}}{\beta_{2m+1}^2}
\cdot
$$
$$
\frac{
\displaystyle
(s_m^2 + 2s_{m-1} s_{m+1} + \sum_{2\le i\le m-1} 2 s_{m-i} s_{m+i})
(s_{m+1}^2 + 2s_m s_{m+2} + \sum_{2\le i\le m} 2 s_{m+1-i} s_{m+1+i})
}
{
\displaystyle 
(2 s_m s_{m+1} + 2 s_{m-1} s_{m+2} + \sum_{2\le i\le m-1} 2 s_{m-i} s_{m+i+1})^2
},
$$
so
\begin{equation}
\label{odd1}
R_{2m+1} = \frac{1}{4}
\frac{ \beta_{2m} \beta_{2m+2}}{\beta_{2m+1}^2}
\cdot 
\end{equation}
$$
\frac{
\begin{array}{c}
\displaystyle
(
1 + 2 R_m + \sum_{2\le i\le m-1} 2 R_{m-(i-1)} \cdots R_m^i \cdots R_{m+(i-1)}
)\cdot \\
\displaystyle
(
1 + 2 R_{m+1} + \sum_{2\le i\le m} 2 R_{m+1-(i-1)} \cdots R_{m+1}^i 
\cdots R_{m+1+(i-1)}
)
\end{array}
}
{
\displaystyle
(1 + R_m R_{m+1} + \sum_{2\le i\le m-1} R_{m-(i-1)}\cdots R_m^i 
R_{m+1}^i \cdots R_{m+i})^2
}.
$$

We now prove by induction that $R_n\le R_{+}=\frac{1}{2}$ for all $n$.
For $2\le n\le 4$ this can be proven by direct computation.
Otherwise, fix some $n\ge 5$, and assume that $R_2$, \dots, $R_{n-1}\le
\frac{1}{2}$.
Set 
$$
\xi_n = \frac{\beta_{n-1} \beta_{n+1}}{\beta_n^2}.
$$ 
It is easy to show that (since $n\ge 5$) $\xi_n\le \frac{10}{19}$.
Now first suppose that $n$ is even, so $n=2m$, $m\ge 3$.
Then from (\ref{even1}), we immediately have
$$
\label{a}
\frac{R_{2m}}{4 \xi_{2m}}
\le
R_m
\frac{
(1+R_{+} R_m + \sum_{i\ge 2} R_{+}^{i(i+1)})^2
}
{
(1+2 R_m)^2
}=A, \qquad \hbox{say.}
$$
We will prove that $A\le \frac{19}{80}$; if this is so, then
$$
R_{2m}\le 4 \xi_{2m} \frac{19}{80} \le 4 \cdot \frac{10}{19} \cdot \frac{19}{80}
= \frac{1}{2},
$$
as desired.
However, 
$\sum_{i\ge 2} R_{+}^{i(i+1)}\le 0.1,$ so
we will have $A\le\frac{19}{80}$ provided that
$$
(\frac{1}{2} R_m+1.1)^2 R_m\le \frac{19}{80} (1+2 R_m)^2,
$$
i.e., if
$$
\frac{1}{4} R_m^3 + 1.1 R_m^2 + 1.21 R_m\le \frac{19}{80} +
\frac{19}{20} R_m + \frac{19}{20} R_m^2.
$$
Since $R_m\le R_{+}= \frac{1}{2}$, it will do to have
$$
(\frac{19}{80} - \frac{1}{4} \cdot (\frac{1}{2})^3) 
+ (\frac{19}{20} - 1.21) R_m + 
(\frac{19}{20} - 1.1) R_m^2\ge 0,
$$
which is true as the expression is positive for $R_m=\frac{1}{2}$ and 
nonincreasing for $R_m\in[0,\frac{1}{2}]$.

The other possibility is that $n=2m+1$ is odd and $m\ge 2$.  Then
from (\ref{odd1}),
$$
\label{b}
\frac{R_{2m+1}}{\frac{1}{4} \xi_{2m+1}}
\le
\frac{
(1+2 R_m + \sum_{i\ge 2} 2 R_{+}^{i^2}) 
(1+2 R_{m+1} + \sum_{i\ge 2} 2 R_{+}^{i^2}) 
}
{
(1+R_m R_{m+1})^2
}=B, \qquad \hbox{say.}
$$
Since $\sum_{i\ge 2} 2 R_{+}^{i^2}\le 0.13$,
we have
\begin{eqnarray*}
B&\le& 
\frac{(1.13+2R_m)(1.13+2R_{m+1})}{1+2R_m R_{m+1}}\\
&=&
\frac{1.13^2+2.26(R_m+R_{m+1})+4R_m R_{m+1}}{1+2 R_m R_{m+1}}\\
&\le& 1.13^2+2.26+(4-2\cdot 1.13^2) \frac{R_m R_{m+1}}{1+2R_m R_{m+1}}\\
&\le& 1.13^2+2.26+(4-2\cdot 1.13^2) \frac{1/4}{1+2/4}\\
&=& 1.13^2+2.26+(4-2\cdot 1.13^2) \frac{1}{6}\le 3.8,
\end{eqnarray*}
so
$$
R_{2m+1}\le 3.8 \frac{1}{4} \xi_{2m+1}\le 3.8 
\cdot \frac{1}{4} \cdot \frac{10}{19}=\frac{1}{2},
$$
as desired.
\end{proof}

\begin{theorem}
Suppose that $q=2$.  Then if $n\ge 2$, $n=qx+y$, and $0\le y\le q$, 
we have
\begin{equation}
\label{eapprox}
\left.{\sum_{\abs{b}=n}}'\prod_{i=0}^{q-1}s_{b_i}\right/
s_x^{q-y} s_{x+1}^y\le 3.
\end{equation}
\end{theorem}
\begin{proof}
Using \eqref{e0}, \eqref{e2}, \eqref{e3} and Lemma \ref{lemr}, we find that
$$
s_{x-i} s_{x+i}\le 2^{-i^2} s_x^2\qquad (x\ge i+1, \ i\ge 0)
$$
and 
$$
s_{x-i} s_{x+i+1}\le 2^{-i(i+1)} s_x s_{x+1} \qquad (x\ge i+1,\ i\ge 0).
$$
It follows that, if $n=2x$ is even,
$$
\sum_{j+k=n,\ j,\,k>0} 
s_j s_k \le s_x^2 (1 + 2 \sum_{i>0} 2^{-i^2}) \le 3 s_x^2,
$$
and if $n=2x+1$ is odd,
$$
\sum_{j+k=n,\ j,\,k>0} 
s_j s_k \le s_x s_{x+1} \cdot 2 \sum_{i\ge 0} 2^{-i(i+1)} \le 3 s_x s_{x+1}.
$$
This completes the proof.
\end{proof}

If $q=2$, we therefore have

\begin{eqnarray}
\label{recursion2}
\log_q s_n&=&\log_q \bt_n
+(q-y)\log_q s_x+y\log_q s_{x+1}+O(1),\\
&&\qquad \hbox{where $n\ge 2$, $n=qx+y$, $0\le y\le q$.}\nonumber
\end{eqnarray}

We conjecture that such a recursion also holds when $q>2$ is integral.
We now show how to reduce recursions like those above to the 
simpler recursion of Section~\ref{recursion}.

\begin{lemma}
\label{sumit}
Fix $q$ and some constant $\bar C$, and for nonnegative integers $n$, set 
$$
\begin{array}{ccll}
w_n&=&-\frac{n^2}{2(1-q^{-1})}-\frac12n\log_q n, & \qquad \hbox{if $n>0$,}\\
w_n&=&0,&\qquad \hbox{if $n=0$.}
\end{array}
$$
Then there is some constant
$C'$ such that if $a_1$, \dots, $a_q$ satisfy $a_1+\dots+a_q=n\ge 1$, 
and for $i=1$, \dots, $q$, we have $|a_i-n/q|<\bar C$, then 
\begin{equation}
\label{bunbun}
|w_n-\log_q \bt_n - (w_{a_1}+\cdots+w_{a_q})|<C'.
\end{equation}
\end{lemma}
\begin{proof}
Write $a_i=n/q+\epsilon_i$.  It will suffice to prove (\ref{bunbun}) for
large $n$.  Take $n$ large enough so that $|\epsilon_i|<\bar C<n/2q$.
Now
\begin{eqnarray*}
w_{a_i}&=&-\frac{(n/q+\epsilon_i)^2}{2(1-q^{-1})}
-\frac12(\frac nq+\epsilon_i)\log_q(\frac nq+\epsilon_i)\\
&=&-\frac{n^2}{2(q^2-q)}-\epsilon_i\frac{n/q}{1-q^{-1}}
-\frac{\epsilon_i^2}{2(1-q^{-1})}\\
&&\qquad -\frac12(\frac nq+\epsilon_i)\log_q \frac nq
-\frac12(\frac nq+\epsilon_i)\log_q(1+\frac{\epsilon_i}{n/q}).
\end{eqnarray*}
Summing over $i$, we get
\begin{eqnarray}
\label{ezoe}
\sum_{1\le i\le q} w_{a_i}
&=&-\frac{n^2}{2(q-1)}-\frac12n(-1+\log_q n)\\
&&\nonumber\qquad
-\sum_{1\le i\le q} \frac{\epsilon_i^2}{2(1-q^{-1})}
+\frac12 (\frac nq+\epsilon_i) \log_q(1+\frac{\epsilon_i}{n/q}).
\end{eqnarray}
By looking at the power series for $\log(1+x)$, we find that
$|\log_q(1+\chi)|\le 4|\chi|$, if $|\chi|<\frac 12$.
By our assumption on $n$, $|\epsilon_i|/(n/q)<\frac12$, so
$|\log_q(1+\epsilon_i/(n/q))|\le 4|\epsilon_i|/(n/q)$.
It follows that the absolute value of the
sum on the right-hand side of (\ref{ezoe})
is bounded, say by $C''$, so
$$
\left|-\frac{n^2}{2(q-1)}+\frac12n-\frac12n\log_q n
-(w_{a_1}+\cdots+w_{a_q})\right|\le C''.
$$
The result now follows from the definition of $w_n$ and
part (i) of Lemma \ref{beta}.
\end{proof}

The significance of Lemma \ref{sumit} is that we can, given a recursion
like (\ref{recursion1}) or (\ref{recursion2}), subtract off $w_n$
from the unknown to obtain a recursion of the form
considered in Section~\ref{recursion}, and then apply Theorem~\ref{qwerty}.
Applying this to \eqref{recursion2} in the $q=2$ case gives us
\eqref{e77}, and concludes the proof of Theorem~\ref{thm:moews}.

If (\ref{recursion2}) were to hold when $q>2$, 
applying Lemma \ref{sumit} to \eqref{recursion2}
and then applying Theorem \ref{qwerty} would give the 

\begin{conjecture}
\label{c26}
If $q>2$ is integral, we have 
\begin{equation}
\label{zzz}
\log_q r_n= -\frac{n^2}{2(q-1)} -\frac12n\log_q n +
(\frac12+ \bar W_q(\log_q n))n + O(1),
\end{equation}
for some continuous $\bar W_q$ of period~1.
\end{conjecture}

Finally, we cannot resist remarking that
applying Lemma \ref{sumit} to (\ref{recursion1})
and then applying Theorem \ref{qwerty} to solve the resulting recursion
of form (\ref{recursion3}) shows that 
$$
\log_q\TTT_n= -\frac{n^2}{2(1-q^{-1})} -\frac12 n\log_q n + 
W_q(\log_q n)n+O(1),
$$
for some continuous $W_q$ of period~1.

This could be used to give an alternate proof of Theorem \ref{thm:moews}.
Since $\TTT_n$ is the greatest tree contribution to $r_n$, and
since the $\log_q$ of the number of labelled $q$-trees with $n$ 
leaves is $O(n)$ (Lemma \ref{blah}), it follows that
$$
\log_q r_n= -\frac{n^2}{2(q-1)} -\frac12n\log_q n +O(n).
$$

We show graphs of $W_2$ and $\bar W_2$ below.

\begin{figure}[p]
\begin{center}
\psfig{file=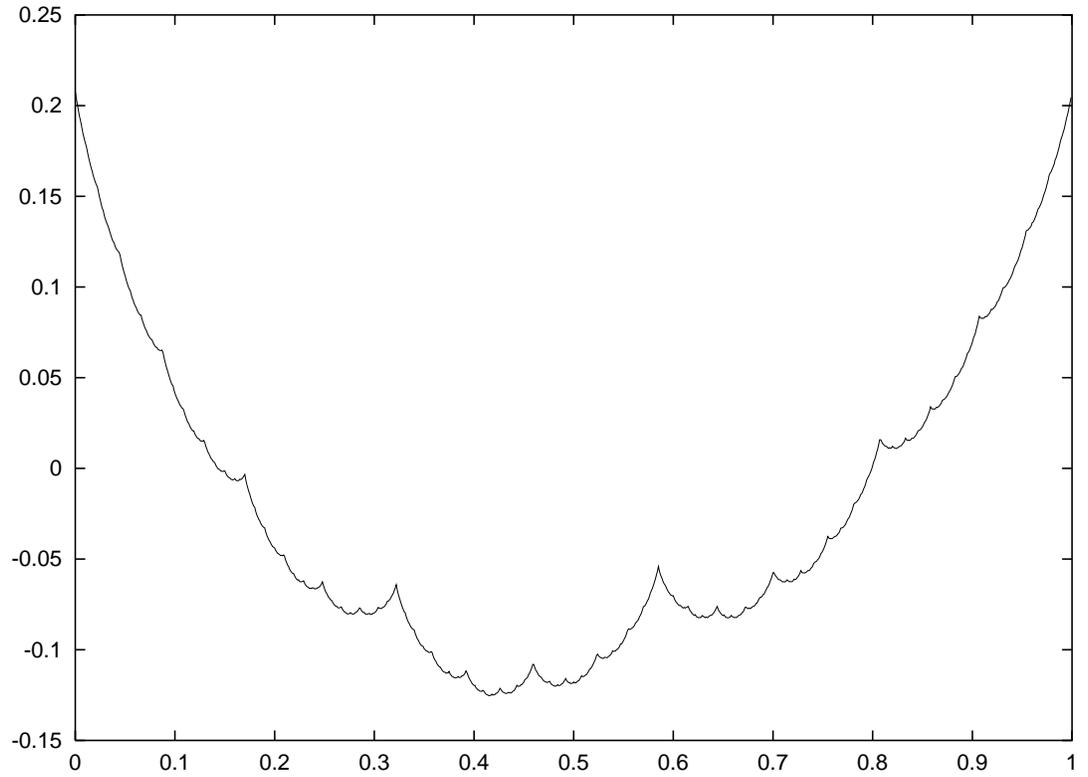,angle=270,width=429bp}
\caption{\label{fig1} Graph of $W_2$ over its period.}
\end{center}
\end{figure}

\begin{figure}[p]
\begin{center}
\psfig{file=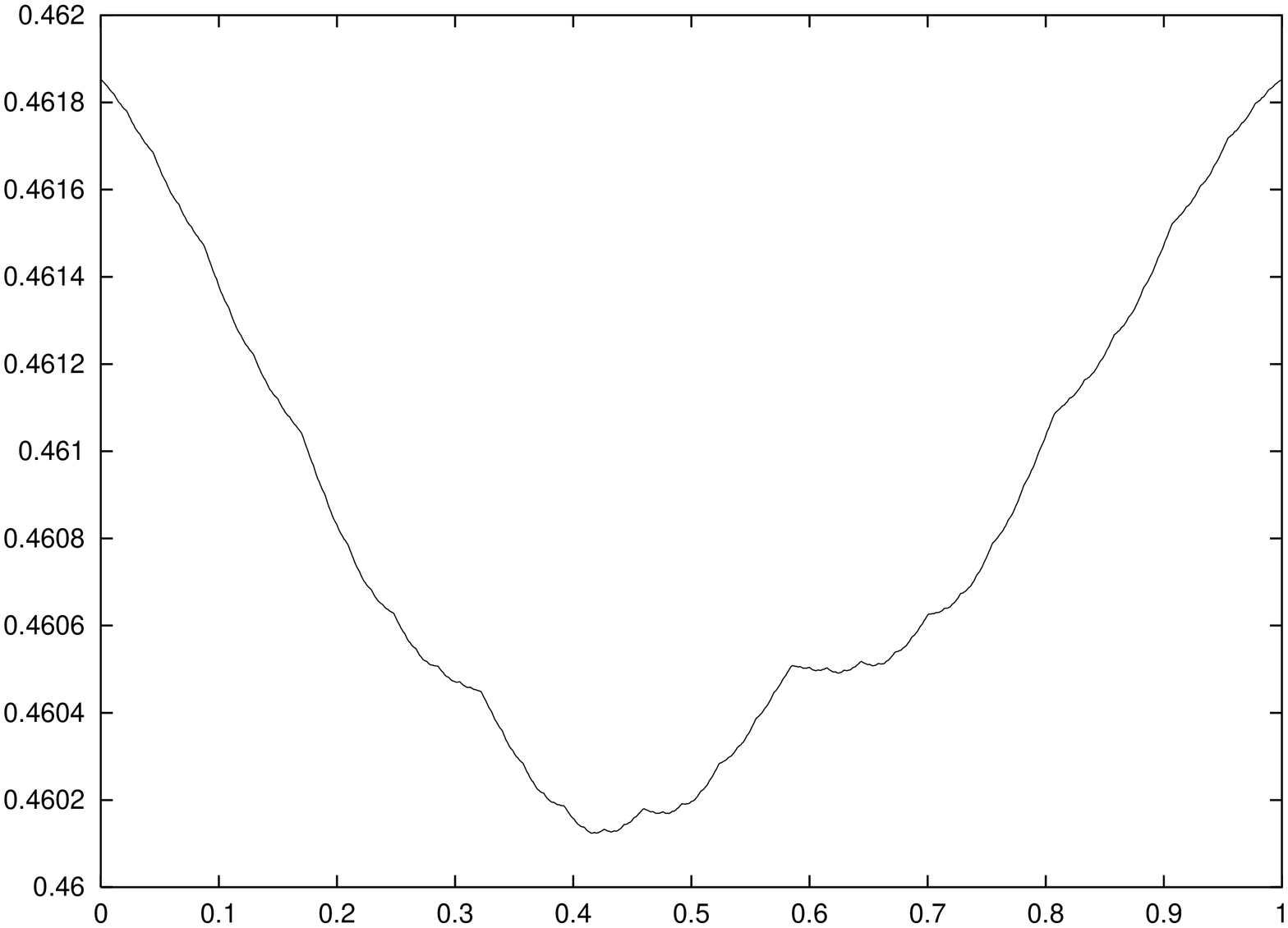,angle=270,width=429bp}
\caption{\label{fig2} Graph of $\bar W_2$ over its period.}
\end{center}
\end{figure}

\end{document}